# Kernel regression uniform rate estimation for censored data under $\alpha$-mixing condition

Zohra Guessoum[1] and Elias Ould-Saïd[2]

[1]*Faculté de Mathématiques, Univ. des Sci. et Tech. Houari Boumédienne, BP 32, El Alia, 16111, Algeria.*
*e-mail:* z0guessoum@hotmail.com

[2]*L.M.P.A. J. Liouville, Univ. du Littoral Côte d'Opale, BP 699, 62228 Calais, France.*
*e-mail:* ouldsaid@lmpa.univ-littoral.fr *(corresponding author)*

**Abstract:** In this paper, we study the behavior of a kernel estimator of the regression function in the right censored model with $\alpha$-mixing data . The uniform strong consistency over a real compact set of the estimate is established along with a rate of convergence. Some simulations are carried out to illustrate the behavior of the estimate with different examples for finite sample sizes.

**Keywords and phrases:** Censored data · Kernel estimator · Nonparametric regression · Rate of convergence · Strong consistency · Strong mixing.
**Mathematics Subject classification (2000)** · 62G05 · 62G20.

## Contents

## 1. Introduction.

Consider a real random variable (rv) $Y$ and a sequence of strictly stationary rv's $(Y_i)_{i\geq 1}$ with common unknown absolutely continuous disribution function (df) $F$. In survival analysis, the rv's may be the lifetimes of patients under study. Let $(C_i)_{i\geq 1}$, be a sequence of censoring rv's with common unknown df $G$. In contrast to statistics for complete data studies, censored model involves pairs $(T_i, \delta_i)_{i=1,\ldots,n}$ where only $T_i = Y_i \wedge C_i$ and $\delta_i = \mathbb{1}_{\{Y_i \leq C_i\}}$ are observed.

Let $X$ be an $\mathbb{R}^d$-valued random vector. Let $(X_i)_{i \geq 1}$ be a sequence of copies of the random vector $X$ and denote by $X_{i,1}, \cdots, X_{i,d}$ the coordinates of $X_i$. The study we perform below is then on the set of observations $(T_i, \delta_i, X_i)_{i \geq 1}$. In regression analysis one expectes to identify, if any, the relationship between the $Y_i$'s and $X_i$'s. This means looking for a function $m^*(X)$ describing this relationship that realizes the minimum of the mean squared error criterion. It is well known that this minimum is achieved by the regression function $m(x)$ defined on $\mathbb{R}^d$ by

$$m(x) = \mathbb{E}(Y \mid X = x).$$

There is a wide range of literature on nonparametric estimation of the regression function and many nonlinear smoothers including kernel, spline, local polynomial, orthogonal methods and so on. For an overview on methods and results for both theoretical and application points of view considering independent or dependent case, we refer the reader to Collomb [13], Silverman [45], Härdle [25], Wahba [48], Wand and Jones [47], Masry and Fan [37], Cai [8] and Cai and Ould-Saïd [9].

In the uncensored case, the behavior of nonparametric estimators built upon mixing sequences is extensively studied. The consistency has been investigated by many authors. Without exhaustivity we quote Robinson [39, 40], Collomb [13], Roussas [42] and Laïb [29]. Some other types of dependence structure have been considered. We refer to Yakowitz [49] for Markov chains, Delecroix [18], Laïb and Ould-Saïd [30] for ergodic processes, Hall and Hart [24] for long-range memory processes, Cai and Roussas [11] for associated random variables. Collomb and Härdle [15] obtained the uniform convergence with rates and some other asymptotic results for a family of kernel robust estimators under a $\varphi$-mixing condition, whereas Gonzalez-Manteiga *et al.* [22] developed a nonparametric test, based on kernel smoothers, to decide whether some covariates could be suppressed in a multidimensional nonparametric regression study. Under the $\alpha$-mixing condition, the uniform strong convergence of the Nadaraya-Watson estimator is treated in Doukhan [19], Bosq [4] and Liebscher [34]. Roussas [42] established the consistency with rate of the regression estimator under some summability requirement. Techniques used in the estimation of nonparametric regression are closely related to density estimation; in this case, kernel estimators have been extensively studied: See for example Roussas [41], Tran [44], Vieu [46] and Liebscher [33, 35]. Cai and Roussas [10] established the strong convergence of the kernel density with rate in the muldimensional case, while Tae and Cox [43] established the same result with a slight difference on the rate. Andrews [1] provides a comprehensive set of results concerning the uniform almost sure convergence. Masry [36] derived sharp rates for the same kind of convergence, but confined attention to the case of bounded regression.

Our goal is to establish the strong uniform convergence with rate for the kernel regression estimate under $\alpha$-mixing condition in random censorship models. For this kind of model, Cai [5, 7] established the asymptotic properties of the Kaplan-Meier estimator. The strong convergence of a hazard rate estimator was examined by Lecoutre and Ould-Saïd [31] while Liebscher [35] derive a rate uniform for the strong convergence of kernel density and hazard rate estimators. His result represents an improvement of that given in Cai [6]. The consistency results concerning the nonparametric estimates of the conditional survival function introduced by Beran [2], Dabrowska [16, 17] in the iid case, were extended by Lecoutre and Ould-Saïd [32] to the strong mixing case. We point out that, for the independent case, the behavior of the regression function under censorship model has been extensively studied. We can quote Carbonez *et al.* [12], Köhler *et al.* [28] and Guessoum and Ould-Saïd [23]. However, few papers deal with the regression function under censoring in the dependent case.

To this end, we were interested in extending the result of Guessoum and Ould-Saïd [23] from the iid to the dependent case. The paper is organized as follows: In Section 2 we give some definitions and notations under the censorship model of the regression function and strong-mixing process. Section 3 is devoted to the assumptions and main result. In Section 4, some simulations are drawn to lend further support to our theoretical results. Proof with auxiliary results are relegated to Section 5.

## 2. Definition of estimators

Suppose that $\{Y_i, i \geq 1\}$ and $\{C_i, i \geq 1\}$ are two independent sequences of stationary random variables. We want to estimate $m(x) = \mathbb{E}(Y | X = x)$ which can be written as $m(x) = \dfrac{r_1(x)}{\ell(x)}$ where

$$r_1(x) = \int_{\mathbb{R}} y f_{X,Y}(x,y) dy \qquad (1)$$

with $f_{\cdot,\cdot}(x,y)$ being the joint density of $(X, Y)$ and $\ell(\cdot)$ the density function of the covariates.

Now, it is well known that the kernel estimator of the regression function

$m(\cdot)$ under censorship model (see, eg Carbonez *et al.* [12]) is given by

$$\tilde{m}_n(x) = \sum_{i=1}^{n} W_{in}(x) \frac{\delta_i T_i}{\bar{G}(T_i)} \qquad (2)$$

where $\bar{G}$ is the survival function of the rv $C$ and

$$W_{in}(x) = \frac{K_d\left(\frac{x-X_i}{h_n}\right)}{\sum_{j=1}^{n} K_d\left(\frac{x-X_j}{h_n}\right)}$$

are the Watson-Nadaraya weights, $K_d$ is a probability density function (pdf) defined on $\mathbb{R}^d$ and $h_n$ a sequence of positive numbers converging to 0 as $n$ goes to infinity. Then (2) can be written

$$\tilde{m}_n(x) =: \frac{\tilde{r}_{1,n}(x)}{\ell_n(x)}$$

where

$$\tilde{r}_{1,n}(x) = \frac{1}{nh_n^d} \sum_{i=1}^{n} \frac{\delta_i T_i}{\bar{G}(T_i)} K_d\left(\frac{x-X_i}{h_n}\right) \quad \text{and} \quad \ell_n(x) = \frac{1}{nh_n^d} \sum_{i=1}^{n} K_d\left(\frac{x-X_i}{h_n}\right). \qquad (3)$$

In practice, $G$ is usually unknown, we replace it by the corresponding Kaplan-Meier [27] estimator (KME) $G_n$ defined by

$$\bar{G}_n(t) = \begin{cases} \prod_{i=1}^{n}\left(1 - \frac{1-\delta_i}{n-i+1}\right)^{\mathbf{1}_{\{Y_i \leq t\}}}, & \text{if } t < Y_{(n)}, \\ 0, & \text{if } t \geq Y_{(n)}. \end{cases}$$

The properties of the KME for dependent variables can be found in Cai [5, 7]. Then a feasible estimator of $m(x)$ is given by: $m_n(x) = \frac{r_{1,n}(x)}{\ell_n(x)}$ where

$$r_{1,n}(x) = \frac{1}{nh_n^d} \sum_{i=1}^{n} \frac{\delta_i T_i}{\bar{G}_n(T_i)} K_d\left(\frac{x-X_i}{h_n}\right) \qquad (4)$$

is an estimator of $r_1(x)$ and $\ell_n(x)$ (defined in (3)) an estimator of $\ell(x)$.
In what follows, we define the endpoints of $F$ and $G$ by $\tau_F = \sup\{y, \bar{F}(y) > 0\}$, $\tau_G = \sup\{y, \bar{G}(y) > 0\}$ and we assume that $\tau_F < \infty$ and $\bar{G}(\tau_F) > 0$ (this implies $\tau_F < \tau_G$).
For technical reasons (see Lemma 5.1), we assume that $\{C_i, i \geq 1\}$ and $\{(X_i, Y_i), i \geq 1\}$ are independent; furthermore this condition is plausible whenever the censoring is independent of the characteristics of the patient under study. We point out that since $Y$ can be a lifetime we can suppose it bounded. We put $\|t\| = \sum_{j=1}^{d} |t_j|$ for $t \in \mathbb{R}^d$.

In order to define the $\alpha$-mixing property, we introduce the following notations. Denote by $\mathcal{F}_i^k(Z)$ the $\sigma$−algebra generated by $\{Z_j, i \leq j \leq k\}$.

**Definition 2.1** *Let* $\{Z_i, i = 1, 2, ...\}$ *denote a sequence of rv's. Given a positive integer $n$, set*

$$\alpha(n) = \sup \left\{ |\mathbb{P}(A \cap B) - \mathbb{P}(A)\mathbb{P}(B)| : A \in \mathcal{F}_1^k(Z) \text{ and } B \in \mathcal{F}_{k+n}^\infty(Z), \ k \in \mathbb{N}^* \right\}.$$

*The sequence is said to be $\alpha$-mixing (strong mixing) if the mixing coefficient $\alpha(n) \to 0$ as $n \to \infty$.*

There exists many processes fulfilling the strong mixing property. We quote, here, the usual ARMA processes which are geometrically strongly mixing, *i.e.*, there exist $\rho \in (0, 1)$ and $a > 0$ such that, for any $n \geq 1$, $\alpha(n) \leq a\rho^n$ (see, *e.g.*, Jones [26]). The threshold models, the EXPAR models (see, Ozaki [38]), the simple ARCH models (see Engle [20]), their GARCH extension (see Bollerslev [3]) and the bilinear Markovian models are geometrically strongly mixing under some general ergodicity conditions.

We suppose that the sequences $\{Y_i, i \geq 1\}$ and $\{C_i, i \geq 1\}$ are $\alpha$-mixing with coefficients $\alpha_1(n)$ and $\alpha_2(n)$, respectively. Cai ([7], Lemma 2) showed that $\{T_i, i \geq 1\}$ is then strongly mixing, with coefficient

$$\alpha(n) = 4 \ \max(\alpha_1(n), \alpha_2(n)).$$

From now on, we suppose that $\{(T_i, \delta_i, X_i) \ i = 1, ..., n\}$ is strongly mixing. Now we are in position to give our assumptions and main result.

## 3. Assumptions and main result

Let $\mathcal{C}$ be a compact set of $\mathbb{R}^d$ which is included in $\mathcal{C}_0 = \left\{ x \in \mathbb{R}^d | \quad \ell(x) > 0 \right\}$. We will make use of the following assumptions gathered here for easy reference:

**A1)** The bandwidth $h_n$ satisfies: $\lim_{n \to +\infty} nh_n^d = +\infty$ and $\lim_{n \to +\infty} h_n^\mu \log \log n = 0$ where $0 < \mu < d$.

**A2)** The kernel $K_d$ is bounded and satisfies:
  i)  $\int_{\mathbb{R}^d} \|t\| K_d(t) dt < +\infty$,
  ii) $\int_{\mathbb{R}^d} (t_1 + t_2 + ... + t_d) K_d^2(t) dt < +\infty$ and $\int_{\mathbb{R}^d} K_d^2(t) dt < +\infty$,
  iii) $\forall (t, s) \in \mathcal{C}^2 \quad |K_d(t) - K_d(s)| \leq \|t - s\|^\gamma$ for $\gamma > 0$.

**A3)** The mixing coefficient $\alpha$ is such that $\alpha(n) = O(n^{-\nu})$ for some $\nu > p + \sqrt{p^2 + 3(d-1)}$ where $p = \frac{\gamma(4+d)+d}{2\gamma}$.

**A4)** The function $r_1(\cdot)$ defined in (1) is continuously differentiable.

**A5)** The function $r_2(x) := \int_{\mathbb{R}^d} y^2 \, f_{X,Y}(x,y) dy$ is continuously differentiable.

**A6)** $\exists D > 0$ such that $\sup_{u,v \in \mathcal{C}} |\ell_{ij}(u,v) - \ell(u)\ell(v)| < D$ where $\ell_{ij}$ is the joint distribution of $(X_i, X_j)$.

**A7)** $\exists \theta > 0, \exists c_1 > 0, \exists c_2 > 0$, such that

$$c_1 n^{\frac{\gamma(3-\nu)}{d[\gamma(\nu+1)+2\gamma+1]} + \theta d} \leq h_n^d \leq c_2 n^{\frac{d}{1-\nu}}.$$

**A8)** The marginal density $\ell(.)$ is continuously differentiable and there exists $\xi > 0$ such that $\ell(x) > \xi \quad \forall x \in \mathcal{C}$.

**Remark 3.1** *Assumption **A1** is very common in functional estimation both in independent and dependent cases. However, it must be reinforced by Assumptions **A3** and **A7** which ensure a practical calculus of the covariance's terms and the convergence of the series which appear in proof of Lemma 3. Assumptions **A2**, **A4**, **A5** and **A6** are needed in the study of the bias term of $r_{1,n}(x)$ which is the kernel estimator of $r_1(x)$. We point out that we do not require for $K_d$ to be symmetric as in Guessoum and Ould-Saïd [23]. Assumption **A8** intervenes in the convergence of the kernel density. Finally, the boundeness of $Y$ is made only for the simplification of the proof. It can be dropped while using truncation methods as in Laïb and Ould-Saïd [30].*

In the sequel letter $C$ denotes any generic constant.
Our main result is given in the following theorem which concerns the rate of the almost sure uniform convergence of the regression function.

**Theorem 3.1** *Under Assumptions **A1-A8**, we have*

$$\sup_{x \in \mathcal{C}} |m_n(x) - m(x)| = O\left(\max\left\{\sqrt{\frac{\log n}{nh_n^d}}, h_n\right\}\right) \quad a.s \quad as \ n \to \infty.$$

**Remark 3.2** *The rate obtained here is slightly different from that obtained by Guessoum and Ould-Saïd [23] in the independent case, which is $O\left(\max\left\{\sqrt{\frac{\log \log n}{nh_n}}, h_n^2\right\}\right)$ for $d = 1$. Their result can easily be generalized for higher dimensional covariate, ie $X \in \mathbb{R}^d$, by adapting their Assumptions **A1** and **A2** to obtain the rate $O\left(\max\left\{\sqrt{\frac{\log \log n}{nh_n^d}}, h_n^2\right\}\right)$.*

**Remark 3.3** *If we choose $h_n = O\left((\log n/n)^{1/d+2}\right)$ then Theorem 1 becomes:*

$$\sup_{x \in \mathcal{C}} |m_n(x) - m(x)| = O\left(\left(\frac{\log n}{n}\right)^{\frac{1}{d+2}}\right) \quad a.s.$$

*This is the optimal rate obtained by Liebscher [34] in the uncensored case.*

## 4. Simulations Study

First, we consider the strong mixing bidimentionnal process generated by:

$$\begin{aligned} X_i &= \rho X_{i-1} + \sqrt{1-\rho^2}\epsilon_i, \\ Y_i &= X_{i+1} \quad i=1,2,...,n, \end{aligned}$$

where $0 < \rho < 1$, $(\epsilon_i)_i$ is a white noise with standard Gaussian distribution and $X_0$ is a standard Gaussian rv independent of $(\epsilon_i)_i$. We also simulate $n$ iid rv $C_i$ exponentially distributed with parameter $\lambda = 1.5$. It is clear that the process $(X_n, Y_n, C_n)$ is stationary and strongly mixing, in fact the process $(X_n)$ is an AR(1) and given $X_1 = x$, we have $Y_1 = \rho x + \sqrt{1-\rho^2}\epsilon_2$, then, $Y_1 \hookrightarrow N(\rho x, 1-\rho^2)$. In all cases we took $\rho = 0.9$. We calculate our estimator based on the observed data $(X_i, T_i, \delta_i,)\ i = 1,...,n$, by choosing a Gaussian kernel $K$. In this case, we have $m(x) = \mathbb{E}(Y_1 | X_1 = x) = \rho x$. In all cases we took $h_n$ satisfying **A1** and **A7**, that is $h_n = O\left((\log n/n)^{1/3}\right)$.

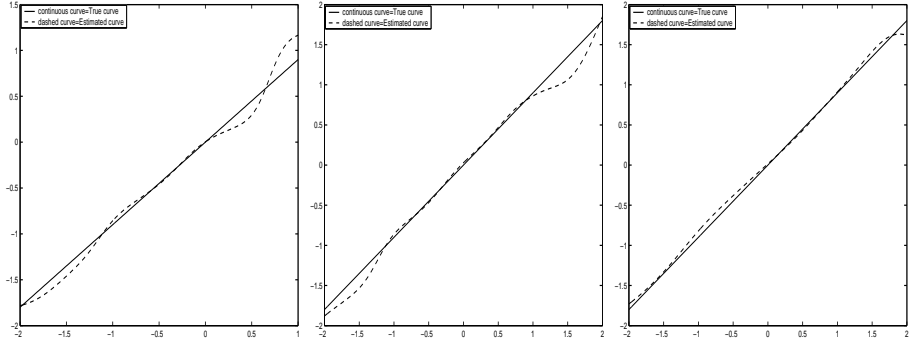

FIG 1. $m(x) = \rho x$, with $n = 50, 100$ and $300$, respectively.

We notice that the quality of fit increases with $n$ (see Figure 1).

We also consider two nonlinear cases

$$Y_i = \sin(\frac{\pi}{2}X_i), \qquad \text{sinus case,} \tag{5}$$

$$Y_i = \frac{5}{12}X_{i+1}^2 - 2, \qquad \text{parabolic case.} \tag{6}$$

Then we have $m(x) = \sin(\frac{\pi}{2}x)$ for (5) and $m(x) = \frac{5}{12}\rho^2 x^2 + \frac{5}{12}(1-\rho^2) - 2$ for (6).

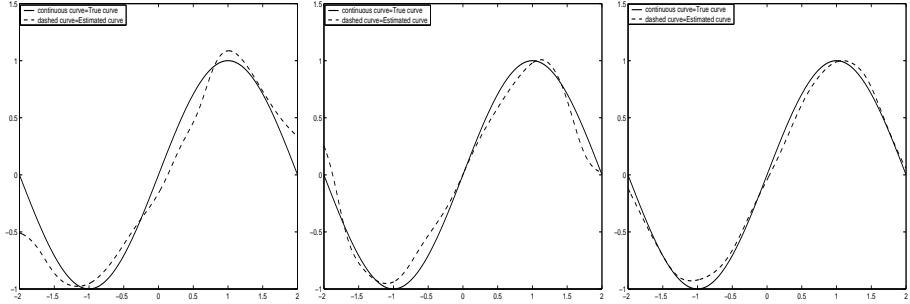

FIG 2. $m(x) = \sin \frac{\pi}{2}x$, with $n = 50$, $100$ and $300$, respectively.

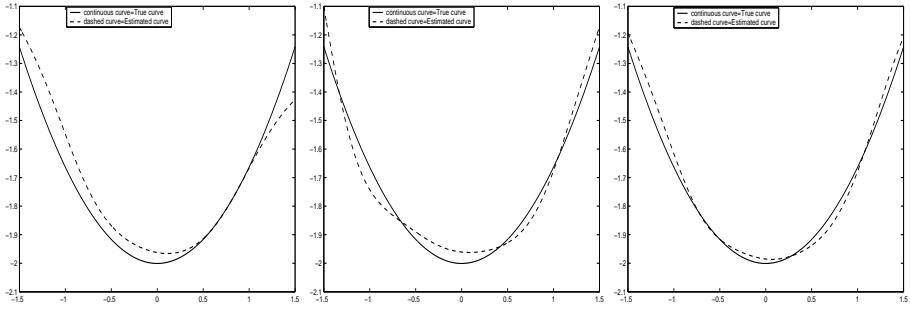

FIG 3. $m(x) = \frac{5}{12}\rho^2 x^2 + \frac{5}{12}(1-\rho^2) - 2$, with $n = 50$, $100$ and $300$, respectively.

Figures 2 and 3 show that the quality of fit for the non linear model is as good as in the linear model.

## 5. Proofs

We split the proof of the Theorem 3.1 in the following Lemmata.

**Lemma 5.1** *Under Assumptions **A1**, **A2** i) and **A4**, for $n$ large enough:*

$$\sup_{x \in \mathcal{C}} |\mathbb{E}\left(\tilde{r}_{1,n}(x)\right) - r_1(x)| = O(h_n) \quad a.s. \quad n \to +\infty.$$

**Proof of Lemma 5.1:** Observe that

$$\mathbb{E}\left(\frac{\delta_1 T_1}{\bar{G}(T_1)} | X_1 = u\right) = \mathbb{E}\left[\mathbb{E}\left[\frac{\mathbb{1}_{\{Y_1 \leq C_1\}} Y_1}{\bar{G}(Y_1)} | Y_1\right] \mid X_1 = u\right]$$
$$= \mathbb{E}\left[\frac{Y_1}{\bar{G}(Y_1)} \mathbb{E}\left[\mathbb{1}_{\{Y_1 \leq C_1\}} | Y_1\right] \mid X_1 = u\right]$$

$$\begin{aligned} &= \mathbb{E}\left[Y_1|X_1 = u\right] \\ &= m(u). \end{aligned}$$

Then, we have from (3)

$$\begin{aligned} \mathbb{E}\left(\tilde{r}_{1,n}(x)\right) - r_1(x) &= \mathbb{E}\left(\frac{1}{h_n^d}\frac{\delta_1 T_1}{\bar{G}(T_1)}K_d\left(\frac{x - X_1}{h_n}\right)\right) - r_1(x) \\ &= \mathbb{E}\left(\frac{1}{h_n^d}K_d\left(\frac{x - X_1}{h_n}\right)\mathbb{E}\left(\frac{\delta_1 T_1}{\bar{G}(T_1)}|X_1\right)\right) - r_1(x) \\ &= \int_{\mathbb{R}^d}\frac{1}{h_n^d}K_d\left(\frac{x - u}{h_n}\right)m(u)\ell(u)du - r_1(x) \\ &= \int_{\mathbb{R}^d} K_d(t)\left[r_1(x - h_n t) - r_1(x)\right]dt \end{aligned}$$

since $r_1 = m\ell$.
A Taylor expansion gives

$$r_1(x - h_n t) - r_1(x) = -h_n(t_1 \frac{\partial r_1}{\partial x_1}(x') + ... + t_d \frac{\partial r_1}{\partial x_d}(x'))$$

where $x'$ is between $x - h_n t$ and $x$. Then

$$\begin{aligned} \sup_{x \in \mathcal{C}}|\mathbb{E}(\tilde{r}_{1,n}(x)) - r_1(x)| &= \sup_{x \in \mathcal{C}}\left|\int_{\mathbb{R}^d} K_d(t)[r_1(x - h_n t) - r_1(x)]dt\right| \\ &\leq h_n \sup_{x \in \mathcal{C}}\int_{\mathbb{R}^d}\left|K_d(t)(t_1\frac{\partial r_1}{\partial x_1}(x') + ... + t_d\frac{\partial r_1}{\partial x_d}(x'))dt\right|. \end{aligned}$$

Then Assumptions **A1, A2** *i*) and **A4** , give the result. ∎

Now, we introduce the following lemma (see Ferraty and Vieu, [21] Proposition A.11 ii), p. 237).

**Lemma 5.2** *Let $\{U_i, i \in \mathbb{N}\}$ be a sequence of real random variables, with strong mixing coefficient $\alpha(n) = O(n^{-\nu})$, $\nu > 1$, such that $\forall n \in \mathbb{N}, \forall i \in \mathbb{N}, 1 \leq i \leq n$ $|U_i| < +\infty$. Then for each $\varepsilon > 0$ and for each $q > 1$ :*

$$\mathbb{P}\left\{\left|\sum_{i=1}^n U_i\right| > \varepsilon\right\} \leq C\left(1 + \frac{\varepsilon^2}{qS_n^2}\right)^{-\frac{q}{2}} + nCq^{-1}\left(\frac{2q}{\varepsilon}\right)^{\nu+1}$$

*where* $S_n^2 = \sum_i\sum_j |cov(U_i, U_j)|$.

**Lemma 5.3** *Under Assumptions **A1-A7**, we have*

$$\sup_{x \in \mathcal{C}}|\tilde{r}_{1,n}(x) - \mathbb{E}\tilde{r}_{1,n}(x))| = O\left(\sqrt{\frac{\log n}{nh_n^d}}\right) \quad a.s \quad as \ n \to \infty.$$

**Proof of Lemma 5.3:**

$\mathcal{C}$ is a compact set, then it admits a covering $\mathcal{S}$ by a finite number $v_n$ of balls $\mathcal{B}_k(x_k^*, a_n^d)$ centered at $x_k^* = (x_{1,k}^*, ... x_{d,k}^*)$, $k \in \{1, ..., v_n\}$. Then for all $x \in \mathcal{C}$ there exists $k \in \{1, ..., v_n\}$ such that $\|x - x_k^*\| \leq a_n^d$, where $a_n$ verifies $a_n^{d\gamma} = h_n^{d(\gamma+\frac{1}{2})} n^{-\frac{d}{2}}$, ($\gamma$ is the same as in Assumption **A2** *iii*)). Since $\mathcal{C}$ is bounded there exists a constant $M > 0$ such that $v_n \leq \frac{M}{a_n^d}$.

Now we set, for $x \in \mathcal{C}$:

$$\Delta_i(x) = \frac{1}{nh_n^d} \frac{\delta_i T_i}{\bar{G}(T_i)} K_d\left(\frac{x - X_i}{h_n}\right) - \mathbb{E}\left(\frac{1}{nh_n^d} \frac{\delta_1 T_1}{\bar{G}(T_1)} K_d\left(\frac{x - X_1}{h_n}\right)\right).$$

It is obvious that

$$\sum_{i=1}^n \Delta_i(x) = \tilde{r}_{1,n}(x) - \mathbb{E}\tilde{r}_{1,n}(x)).$$

Writing $\Delta_i(x) - \Delta_i(x_k^*) = \tilde{\Delta}_i(x)$, we have clearly $|\Delta_i(x)| \leq |\tilde{\Delta}_i(x)| + |\Delta_i(x_k^*)|$. Then,

$$\sup_{x \in \mathcal{C}} \left|\sum_{i=1}^n \tilde{\Delta}_i(x)\right| \leq \sup_{x \in \mathcal{C}} \left\{\frac{1}{n} \sum_{i=1}^n \frac{\delta_i |T_i|}{\bar{G}(T_i)} \frac{1}{h_n^d} \left|K_d\left(\frac{x - X_i}{h_n}\right) - K_d\left(\frac{x_k^* - X_i}{h_n}\right)\right|\right\}$$
$$+ \sup_{x \in \mathcal{C}} \left\{\mathbb{E}\left(\frac{\delta_1 |T_1|}{\bar{G}(T_1)} \frac{1}{h_n^d} \left|K_d\left(\frac{x - X_1}{h_n}\right) - K_d\left(\frac{x_k^* - X_1}{h_n}\right)\right|\right)\right\}.$$

From Assumption **A2** *iii*)

$$\sup_{x \in \mathcal{C}} \left|\sum_{i=1}^n \tilde{\Delta}_i(x)\right| \leq \sup_{x \in \mathcal{C}} \left(\frac{2\mathbb{E}(|Y_1|)}{\bar{G}(\tau_F)} \frac{1}{h_n^d} \left\|\frac{x - x_k^*}{h_n}\right\|^\gamma\right) \leq \frac{\mathbb{E}(|Y_1|)}{\bar{G}(\tau_F)} \frac{a_n^{d\gamma}}{h_n^{\gamma+d}}$$
$$\leq \frac{\mathbb{E}(|Y_1|)}{\bar{G}(\tau_F)} \frac{h_n^{d(\gamma+\frac{1}{2})} n^{-\frac{d}{2}}}{h_n^{\gamma+d}} \leq \frac{C}{\sqrt{nh_n^d}} h_n^{\gamma(d-1)}.$$

Assumption **A1** implies that $\sup_{x \in \mathcal{C}} \left|\sum_{i=1}^n \tilde{\Delta}_i(x)\right| = O\left(\frac{1}{\sqrt{nh_n^d}}\right)$ a.s.

On the other hand, let $U_i = nh_n^d \Delta_i(x_k^*)$. In order to apply Lemma 5.2, we have to calculate $S_n^2$. It is clear that

$$S_n^2 = \sum_i \sum_{j, i \neq j} |cov(U_i, U_j)| + nVar(U_1).$$

We have
$$Var(U_1) = \mathbb{E}\left[\frac{\delta_1^2 T_1^2}{\bar{G}^2(T_1)} K_d^2\left(\frac{x_k^* - X_1}{h_n}\right)\right] - \mathbb{E}^2\left[\frac{\delta_1 T_1}{\bar{G}(T_1)} K_d\left(\frac{x_k^* - X_1}{h_n}\right)\right] =: I_1 - I_2.$$

Using the conditional expectation properties and a change of variables, we get

$$\begin{aligned} I_1 &= \mathbb{E}\left[\frac{\delta_1^2 T_1^2}{\bar{G}^2(T_1)} K_d^2\left(\frac{x_k^* - X_1}{h_n}\right)\right] \\ &= \mathbb{E}\left[K_d^2\left(\frac{x_k^* - X_1}{h_n}\right) \mathbb{E}\left(\frac{\delta_1^2 T_1^2}{\bar{G}^2(T_1)}\bigg| X_1\right)\right] \\ &\leq \frac{h_n^d}{\bar{G}(\tau_F)} \int_{\mathbb{R}^d} K_d^2(t) r_2(x_k^* - h_n t) dt. \end{aligned}$$

By a Taylor expansion around $x_k^*$, under Assumptions **A2** *ii)* and **A5**, we obtain
$$I_1 = O(h_n^d).$$

From Assumption **A4**,
$$I_2 = \mathbb{E}^2\left[K_d\left(\frac{x_k^* - X_1}{h_n}\right) \mathbb{E}\left[\frac{\delta_1 T_1}{\bar{G}(T_1)}\bigg| X_1\right]\right] = \left[\int_{\mathbb{R}^d} K_d\left(\frac{x_k^* - u}{h_n}\right) r_1(u) dt\right]^2 = O(h_n^{2d}).$$

Finally $Var(U_1) = O(h_n^d)$.
Now let $S_n^{2*} = \sum_i \sum_{j,\ i\neq j} |cov(U_i, U_j)|$, a direct calculus of $|cov(U_i, U_j)|$ gives

$$\begin{aligned} |cov(U_i, U_j)| &= |\mathbb{E} U_i U_j| = \bigg|\mathbb{E}\bigg\{\bigg[\frac{\delta_i T_i}{\bar{G}(T_i)} K_d\left(\frac{x_k^* - X_i}{h_n}\right) - \mathbb{E}\left(\frac{\delta_i T_i}{\bar{G}(T_i)} K_d\left(\frac{x_k^* - X_i}{h_n}\right)\right)\bigg] \\ &\quad \times \bigg[\frac{\delta_j T_j}{\bar{G}(T_j)} K_d\left(\frac{x_k^* - X_j}{h_n}\right) - \mathbb{E}\left(\frac{\delta_j T_j}{\bar{G}(T_j)} K_d\left(\frac{x_k^* - X_j}{h_n}\right)\right)\bigg]\bigg\}\bigg| \\ &\leq \bigg|\mathbb{E}\left(\frac{Y_i}{\bar{G}(Y_i)} K_d\left(\frac{x_k^* - X_i}{h_n}\right) \frac{Y_j}{\bar{G}(Y_j)} K_d\left(\frac{x_k^* - X_j}{h_n}\right)\right) \\ &\quad - \mathbb{E}\left(\frac{Y_i}{\bar{G}(Y_i)} K_d\left(\frac{x_k^* - X_i}{h_n}\right)\right) \mathbb{E}\left(\frac{Y_j}{\bar{G}(Y_j)} K_d\left(\frac{x_k^* - X_j}{h_n}\right)\right)\bigg| \\ &\leq C h_n^{2d} \int_{\mathbb{R}^d}\int_{\mathbb{R}^d} K_d(z) K_d(t) \left[\ell_{ij}\left(x_k^* - zh_n, x_k^* - th_n\right) - \ell_i(x_k^* - zh_n)\ell_j(x_k^* - th_n)\right] dz\ dt. \end{aligned}$$

Assumption **A6** gives
$$|cov(U_i, U_j)| = O(h_n^{2d}). \tag{7}$$

On the other hand, from a result in Bosq ([4], p. 22), we have

$$|cov(U_i, U_j)| \leq C\alpha(|i-j|). \tag{8}$$

Then to evaluate $S_n^{2*}$ the idea is to introduce a sequence of integers $w_n$ which we precise below. Then we use (7) for the close $i$ and $j$ and (8) otherwise. That is

$$\begin{aligned}
S_n^{2*} &= \sum\sum\nolimits_{0<|i-j|\leq w_n} |cov(U_i, U_j)| + \sum\sum\nolimits_{|i-j|>w_n} |cov(U_i, U_j)| \\
&\leq C\sum\sum\nolimits_{|i-j|>w_n}\sum\sum\nolimits_{0<|i-j|\leq w_n} h_n^{2d} + C\alpha(|i-j|) \\
&\leq C\left(nh_n^{2d}w_n\right) + Cn^2\alpha(w_n).
\end{aligned}$$

Now choosing $w_n = \left[\frac{1}{h_n^d}\right] + 1$, we have $S_n^{2*} \leq O(nh_n^d) + Cn^2\alpha\left(\frac{1}{h_n^d}\right)$. Assumption **A3** and the right part of Assumption **A7** yield $n^2\alpha(\frac{1}{h_n^d}) = O(nh_n^d)$. So

$$S_n^{2*} = O(nh_n^d).$$

Finally, we have

$$S_n^2 = S_n^{2*} + nVar(U_1) = O(nh_n^d).$$

Then, for $\varepsilon > 0$, applying Lemma 5.2, we have

$$\mathbb{P}\left\{\left|\sum_{i=1}^n \Delta_i(x^{*k})\right| > \varepsilon\right\} = \mathbb{P}\left\{\left|\sum_{i=1}^n U_i\right| > nh_n^d\varepsilon\right\}$$
$$\leq C\left(1 + C\frac{\varepsilon^2 nh_n^d}{q}\right)^{-\frac{q}{2}} + nCq^{-1}\left(\frac{q}{\varepsilon nh_n^d}\right)^{\nu+1} \tag{9}$$

If we replace $\varepsilon$ by $\varepsilon_0\sqrt{\frac{\log n}{nh_n^d}}$ for all $\varepsilon_0 > 0$ in (9), we get

$$\mathbb{P}\left\{\left|\sum_{i=1}^n \Delta_i(x^{*k})\right| > \varepsilon_0\sqrt{\frac{\log n}{nh_n^d}}\right\} \leq C\left(1 + C\frac{\varepsilon_0^2 \log n}{q}\right)^{-\frac{q}{2}} + nCq^{-1}\left(\frac{q}{\varepsilon_0\sqrt{nh_n^d \log n}}\right)^{\nu+1}. \tag{10}$$

By choosing $q = (\log n)^{1+b}$ $(b > 0)$, (10) becomes

$$\mathbb{P}\left\{\left|\sum_{i=1}^n \Delta_i(x^{*k})\right| > \varepsilon_0\sqrt{\frac{\log n}{nh_n^d}}\right\} \leq Cn^{-C\varepsilon_0^2} + nCq^{-1}\left(\frac{q}{\varepsilon_0}\right)^{\nu+1}\left(nh_n^d \log n\right)^{-\frac{\nu+1}{2}}$$
$$\leq Cn^{-C\varepsilon_0^2} + C\varepsilon_0^{-(\nu+1)}(\log n)^{\nu(1+b)} n^{1-\frac{\nu+1}{2}} h_n^{-\frac{d(\nu+1)}{2}}.$$

Now we can write

$$\mathbb{P}\left\{\max_{k=1,\ldots,v_n}\left|\sum_{i=1}^n\Delta_i(x^{*k})\right|>\varepsilon_0\sqrt{\frac{\log n}{nh_n^d}}\right\} \leq \sum_{i=1}^{v_n}\mathbb{P}\left\{\left|\sum_{i=1}^n\Delta_i(x^{*k})\right|>\varepsilon_0\sqrt{\frac{\log n}{nh_n^d}}\right\}.$$

$$\leq Ma_n^{-d}\left(Cn^{-C\varepsilon_0^2}+C\varepsilon_0^{-(\nu+1)}(\log n)^{\nu(1+b)}n^{1-\frac{\nu+1}{2}}h_n^{-\frac{d(\nu+1)}{2}}\right)$$

$$\leq CMh_n^{-d(1+\frac{1}{2\gamma})}n^{\frac{d}{2\gamma}-C\varepsilon_0^2}$$

$$+ MC\varepsilon_0^{-(\nu+1)}(\log n)^{\nu(1+b)}n^{1-\frac{\nu+1}{2}+\frac{d}{2\gamma}}h_n^{-\frac{d(\nu+1)}{2}-\frac{d}{2\gamma}-d}$$

$$=: CMJ_1 + MC\varepsilon_0^{-(\nu+1)}J_2. \qquad (11)$$

We have from the left part of assumption **A7**

$$J_2 \leq C(\log n)^{\nu(1+b)}n^{1-\frac{\nu+1}{2}+\frac{d}{2\gamma}}n^{-\frac{(3-\nu)}{2}-\theta d(\frac{\gamma(\nu+1)+2\gamma+1}{2\gamma})}$$

$$\leq C(\log n)^{\nu(1+b)}n^{-1-\theta d(\frac{\gamma(\nu+1)+2\gamma+1-\frac{1}{\theta}}{2\gamma})}.$$

Then, for an appropriate choice of $\theta$, $J_2$ is the general term of a convergent series. In the same way, $J_1 \leq n^{\varsigma-C\varepsilon_0^2}$ and we can choose $\varepsilon_0$ such that $J_1$ is the general term of convergent series. Finally, applying Borel-Cantelli lemma, to (11) gives the result. ∎

**Remark 5.1** *We point out that the parameter $\theta$ of Assumption **A7** can be chosen such as:*

$$\frac{1}{\gamma(\nu+1)+2\gamma+1}<\theta<\frac{1}{1-\nu}-\frac{\gamma(3-\nu)}{d[\gamma(\nu+1)+2\gamma+1]}.$$

*This condition ensures the convergence of the series of Lemma 3.*

**Lemma 5.4** *Under Assumptions **A1-A3** and **A6-A8**,*

$$\sup_{x\in\mathcal{C}}|\ell_n(x)-\ell(x)|=O\left(\max\left\{\sqrt{\frac{\log n}{nh_n^d}},h_n\right\}\right)\quad a.s\quad as\ n\to\infty.$$

**Proof of Lemma 5.4:** We have

$$\sup_{x\in\mathcal{C}}|\ell(x)-\ell_n(x)|\leq\sup_{x\in\mathcal{C}}|\ell_n(x)-\mathbb{E}(\ell_n(x))|+\sup_{x\in\mathcal{C}}|\mathbb{E}(\ell_n(x))-\ell(x)|.$$

By Assumptions **A1-A3**, **A6-A8**, by an analogous proof to that of Lemma 5.3 without censoring ( that is $\bar{G}(T_i)=1$, $\delta_i=1$ and $Y_i=1$) and putting

$\varepsilon = \varepsilon_0 \sqrt{\frac{\log n}{nh_n^d}}$ we get

$$\sup_{x \in \mathcal{C}} |\ell_n(x) - \mathbb{E}(\ell_n(x))| = O\left(\sqrt{\frac{\log n}{nh_n^d}}\right). \quad (12)$$

Furthermore, under **A2** *i*) and **A8** and using a Taylor expansion, we get

$$\sup_{x \in \mathcal{C}} |\mathbb{E}(\ell_n(x)) - \ell(x)| = O(h_n)$$

which permit us to conclude. ∎

**Lemma 5.5** *Under Assumptions **A1-A3, A6-A8**, we have*

$$\sup_{x \in \mathcal{C}} |r_{1,n}(x) - \tilde{r}_{1,n}(x)| = o\left(\frac{1}{\sqrt{nh_n^\mu}}\right) \quad a.s \quad as \ n \to \infty.$$

**Proof of Lemma 5.5:** We have from (3) and (4)

$$\begin{aligned}
|r_{1,n}(x) - \tilde{r}_{1,n}(x)| &= \frac{1}{nh_n^d} \left| \sum_{i=1}^n \frac{\mathbb{I}_{\{Y_i < C_i\}} Y_i}{\bar{G}_n(Y_i)} K_d\left(\frac{x - X_i}{h_n}\right) - \frac{\mathbb{I}_{\{Y_1 < C_1\}} Y_i}{\bar{G}(Y_i)} K_d\left(\frac{x - X_i}{h_n}\right) \right| \\
&\leq \frac{1}{nh_n^d} \left| \sum_{i=1}^n Y_i K_d\left(\frac{x - X_i}{h_n}\right) \frac{\bar{G}(Y_i) - \bar{G}_n(Y_i)}{\bar{G}_n(Y_i)\bar{G}(Y_i)} \right| \\
&\leq \frac{1}{\bar{G}_n(\tau_F)\bar{G}(\tau_F)} \sup_{t \leq \tau_F} (|\bar{G}_n(t) - \bar{G}(t)|) \frac{1}{nh_n^d} \sum_{i=1}^n |Y_i| K_d\left(\frac{x - X_i}{h_n}\right).
\end{aligned}$$

In the same way as for Theorem 2 of Cai (2001), it can be shown under **A3** that

$$\sup_{t \leq \tau_F} (|\bar{G}_n(t) - \bar{G}(t)|) = O\left(\sqrt{\frac{\log \log n}{n}}\right) \ a.s.$$ Furthermore, from the definition of $\ell_n(x)$, Lemma 5.4, Assumptions **A1, A2** and **A8**, and the fact that $Y$ is bounded we get the result. ∎

**Proof of Theorem 3.1:** We have

$$\begin{aligned}
\sup_{x \in \mathcal{C}} |m_n(x) - m(x)| &\leq \sup_{x \in \mathcal{C}} \left\{ \left| \frac{r_{1,n}(x)}{\ell_n(x)} - \frac{\tilde{r}_{1,n}(x)}{\ell_n(x)} \right| + \left| \frac{\tilde{r}_{1,n}(x)}{\ell_n(x)} - \frac{\mathbb{E}\tilde{r}_{1,n}(x)}{\ell_n(x)} \right| \right. \\
&\quad + \left. \left| \frac{\mathbb{E}(\tilde{r}_{1,n}(x))}{\ell_n(x)} - \frac{r_1(x)}{\ell_n(x)} \right| + \left| \frac{r_1(x)}{\ell_n(x)} - \frac{r_1(x)}{\ell(x)} \right| \right\} \\
&\leq \frac{1}{\inf \ell_n(x)} \left\{ \sup_{x \in \mathcal{C}} |r_{1,n}(x) - \tilde{r}_{1,n}(x)| + \sup_{x \in \mathcal{C}} |\tilde{r}_{1,n}(x) - \mathbb{E}\tilde{r}_{1,n}(x)| \right. \\
&\quad + \left. \sup_{x \in \mathcal{C}} |\mathbb{E}(\tilde{r}_{1,n}(x) - r_1(x))| + \sup_{x \in \mathcal{C}} \left(|r_1(x)|\xi^{-1}\right) \sup_{x \in \mathcal{C}} |\ell(x) - \ell_n(x)| \right\}. \ (13)
\end{aligned}$$

The kernel estimator $\ell_n(x)$ is almost surely bounded away from 0 because of Lemma 5.4 and the fact that the second part of Assumption A8.

Then, (13) in conjunction with Lemmas 5.1, 5.3, 5.4 and 5.5 we conclude the proof. ∎